\documentclass[1p,times,final]{elsarticle}

\biboptions{semicolon}

\usepackage{framed,multirow}

\usepackage{amssymb}
\usepackage{url}
\usepackage{xcolor}
\definecolor{newcolor}{rgb}{.8,.349,.1}

\usepackage{amsmath}
\usepackage{bm}

\newcommand{\ABS}[1]{\ensuremath{\lvert {#1} \rvert}}
\newcommand{\dpone}[2]{\ensuremath{\displaystyle\frac{\partial {#1}}{\partial {#2}}}}

\newcommand{\bF}{\ensuremath{\bm{F}}}
\newcommand{\bG}{\ensuremath{\bm{G}}}

\newcommand{\bI}{\ensuremath{\bm{I}}}
\newcommand{\bJ}{\ensuremath{\bm{J}}}
\newcommand{\bK}{\ensuremath{\bm{K}}}

\newcommand{\bQ}{\ensuremath{\bm{Q}}}

\newcommand{\bff}{\ensuremath{\bm{f}}}
\newcommand{\bg}{\ensuremath{\bm{g}}}

\newcommand{\bq}{\ensuremath{\bm{q}}}

\newcommand{\bs}{\ensuremath{\bm{s}}}

\newcommand{\bu}{\ensuremath{\bm{u}}}

\newcommand{\bx}{\ensuremath{\bm{x}}}
\newcommand{\by}{\ensuremath{\bm{y}}}

\newcommand{\Fcal}{\ensuremath{\mathcal{F}}}
\newcommand{\Ncal}{\ensuremath{\mathcal{N}}}

\newcommand{\Ocal}{\ensuremath{{\mathcal{O}}}}

\newcommand{\defeq}{\mathrel{\overset{\makebox[0pt]{\mbox{\normalfont\tiny\sffamily def}}}{=}}}

\journal{Journal of Computational Physics}

\begin{document}

	\title{Universal image systems for non-periodic and periodic Stokes flows above a no-slip wall}

	\author[add1,add2]{Wen Yan\corref{cor1}}
	\ead{wyan@flatironinstitute.org}
	\ead{wenyan4work@gmail.com}
	\author[add1,add2]{Michael Shelley}
	\ead{mshelley@flatironinstitute.org}
	%\fntext[flabel1]{Corresponding author}
	%\fntext[flabel2]{textforlabel2}
	\cortext[cor1]{Corresponding author}
	%\cortext[cor2]{textforcor2}
	
	\address[add1]{Center for Computational Biology, Flatiron Institute, Simons Foundation, New York 10010}
	\address[add2]{Courant Institute of Mathematical Sciences, New York University, New York 10010}

\begin{abstract}
		
		It is well-known that by placing judiciously chosen image point forces and doublets to the Stokeslet above a flat wall, the no-slip boundary condition can be conveniently imposed on the wall [Blake, J. R. Math. Proc. Camb. Philos. Soc. 70(2), 1971: 303.].
		However, to further impose periodic boundary conditions on directions parallel to the wall usually involves tedious derivations because single or double periodicity in Stokes flow may require the periodic unit to have no net force, which is not satisfied by the well-known image system.
		In this work we present a force-neutral image system. This neutrality allows us to represent the Stokes image system in a universal formulation for non-periodic, singly periodic and doubly periodic geometries.
		This formulation enables the black-box style usage of fast kernel summation methods.
		We demonstrate the efficiency and accuracy of this new image method with the periodic kernel independent fast multipole method in both non-periodic and doubly periodic geometries.
		We then extend this new image system to other widely used Stokes fundamental solutions, including the Laplacian of the Stokeslet and the Rotne-Prager-Yamakawa tensor.
	
\end{abstract}


\maketitle

\section{Introduction}
\label{sec:intro}
No-slip boundaries in Stokes flow are central to much flow phenomena.
For example, for Brownian suspensions above a no-slip wall, the wall not only constrains the motion of particles, but fundamentally changes the self and collective Brownian motion of suspensions by inducing anisotropy and screening effects in the mobility of particles \citep{Lele_Swan_Brady_Wagner_Furst_2011,Michailidou_Swan_Brady_Petekidis_2013,Usabiaga_Delmotte_Donev_2017}. 
Another example is that swimming microorganisms may swim upstream near a no-slip boundary in an imposed flow due to either hydrodynamic or non-hydrodynamic causes \citep{Kaya_Koser_2012,Spagnolie_Lauga_2012,Wioland_Lushi_Goldstein_2016,Ezhilan_Saintillan_2015}.

To compute the Stokes flow above a no-slip wall, the image method of \citet{Blake_1971} is a popular choice. For a Stokeslet above a wall Blake showed that the no-slip condition was satisfied by adding an image Stokeslet, a modified source doublet, and a modified force doublet to the original Stokeslet. Similar methods have also been developed by \citet{Mitchell_Spagnolie_2015,Mitchell_Spagnolie_2017}. Recently, \citet{Gimbutas_Greengard_Veerapaneni_2015} developed a simpler image system.
This system invokes standard Stokes and Laplace kernel evaluations only, which is compatible with the Fast Multipole Method (FMM).
However, to further impose periodic boundary conditions on the two directions parallel to the no-slip wall is no simple task, because different kernel summations in the image system need to be periodized simultaneously and coupled to each other.
\citet{Nguyen_Leiderman_2015} recently derived the Ewald summation formulation for the doubly periodic Stokeslet image system, but their method show a non-optimal $\Ocal(N^2)$ scaling for $N$ point forces. Their method was recently applied in the study of ciliary beating \citep{Hoffmann_Cortez_2017}. To our knowledge, the singly periodic Stokeslet image system above a no-slip wall has not yet been derived.

For convenience and efficiency, it is desirable to develop an image system where each kernel sum $\bg^t = \sum_s \bK(\bx^t,\by^s)\bq(\by^s)$ can be independently computed and periodized.
Here $\bx^t$ and $\by^s$ are target and source points with indices $s$ and $t$. $\bK$ is the kernel function. The kernel sum could be simply written as $\bg = \bK\bq$, where the indices $s,t$ are suppressed. For Stokes and Laplace kernel sums, recently developed optimal fast periodic kernel summation methods with flexible periodic boundary conditions can be used, including the Spectral Ewald methods \citep{Lindbo_Tornberg_2011,Lindbo_Tornberg_2012,Shamshirgar_Tornberg_2017}, which scale as $\Ocal(N\log N)$, and the periodic Kernel Independent Fast Multipole Method (KIFMM) method by \citet{PBCFMM2018}, which scales as $\Ocal(N)$. 
However the image systems developed by \citet{Gimbutas_Greengard_Veerapaneni_2015} does not work in this framework, because the partially periodic (i.e., simply or doubly periodic) summations for the Stokeslet and the Laplace monopole kernel do not allow a net force or a net monopole in a periodic box, as otherwise the infinite periodic summations diverge. 
Unfortunately this requirement is not satisfied by the image system of \cite{Gimbutas_Greengard_Veerapaneni_2015}.

In this work we propose a new image system for the Stokeslet, which satisfies the neutrality condition by rearranging the Stokeslet and Laplace kernel sums in the image system by \citet{Gimbutas_Greengard_Veerapaneni_2015}. 
Therefore any singly or doubly periodic kernel summation method can be used as a black-box routine to periodize this new image system.

In Section~\ref{sec:formulation} we briefly derive the new image system.
Numerical results for Stokeslet above a no-slip wall with non-periodic and doubly periodic boundary conditions are presented in Section~\ref{sec:results}. 
In Section~\ref{sec:extension} we extend the new image to the Laplacian of Stokeslet and the widely used Rotne-Prager-Yamakawa tensor \citep{Rotne_Prager_1969,Yamakawa_1970}.
We conclude this work with a brief discussion about its coupling to fast summation methods, and its extension to other kernels.

\section{Formulation}
\label{sec:formulation}
We first consider a point force $\bff=(f_1,f_2,f_3)$ located at $\by=(y_1,y_2,y_3)$ above an infinite no-slip wall at the plane $x_3=0$.
We define the image force $\bff^I = (f_1,f_2,-f_3)$ located at $\by^I=(y_1,y_2,-y_3)$ below the wall. 
The complete image system to satisfy the no-slip condition on the wall is given by \citet{Gimbutas_Greengard_Veerapaneni_2015} following the Papkovich-Neuber technique:
\begin{subequations}\label{eq:imageGVG}
\begin{align}
\bu(\bx) &= \bJ(\bx,\by)\bff + \bJ(\bx,\by^I)\left(-\bff^I\right) - \bu^C(\bx),\\
\bu^C(\bx) &= x_3 \nabla_x\phi(\bx) - \hat{\bx}_3\phi(\bx),\quad \phi(\bx)\defeq G^S(\bx,\by^I) f_3^I + \bG^D(\bx,\by^I) (y_3\bff^I ),
\end{align}
\end{subequations}
where $\hat{\bx}_3$ is the unit vector in the $x_3$ direction.
In this expression three kernels are involved: the Laplace monopole kernel 
 $G^S(\bx,\by) = \frac{1}{4\pi}\frac{1}{\ABS{\bx-\by}}$, the Laplace dipole kernel
 $\bG^D(\bx,\by)= \frac{1}{4\pi}\frac{\bx-\by}{\ABS{\bx-\by}^3} = \nabla_y G^S(\bx,\by)$, and the Stokeslet $\bJ(\bx,\by) = \frac{1}{8\pi}\left( \frac{\bI}{\ABS{\bx-\by}} + \frac{(\bx-\by)(\bx-\by)}{\ABS{\bx-\by}^3} \right)$.
We set the fluid viscosity to $\eta=1$ in $\bJ$ for simplicity.
It is clear that the net force is $\bff+(-\bff^I)=(0,0,2f_3)\neq 0$ in the Stokes kernel sum, and the net monopole is $f_3^I = -f_3\neq0$ in the Laplace monopole kernel sum.
This forbids us to apply partially periodic kernel sum methods directly. The requirement of neutrality is straightforward to understand for Laplace kernels. For Stokeslet this depends on the particular periodic boundary condition. With triply periodic periodic boundary condition, the net force within a periodic box does not have to be zero because the net force could be balanced by the global pressure gradient \cite{Hasimoto_1959}. However the net force must be zero with singly and doubly periodic boundary conditions, as demonstrated by \citet{Lindbo_Tornberg_2011}. 

To remove the net force and net monopole, we convert the third component of Stokes force into a Laplace monopole kernel sum following the idea of \citet{Tornberg_Greengard_2008}. 
This involves tedious algebraic manipulations and we only summarize the results here.
The new image system splits the flow velocity into $4$ independent parts $\bu(\bx) = \bu^S + \bu^D +\bu^{L1} +\bu^{L2}$, where each part is computed by one kernel sum. In the following, $\bff_{xy}=(f_1,f_2,0)$ denotes the $x_1,x_2$ component of the point force $\bff$, parallel to the no-slip wall.
\begin{subequations}\label{eq:imageNew}
\begin{align}
\bu^S&=\bJ(\bx,\by)\bff_{xy} + \bJ(\bx,\by^I)\left(-\bff_{xy}\right),\\
\bu^D& = \left(x_3 \nabla_x - \hat{\bx}_3\right)\phi^D(\bx),\quad\text{ with } \phi^D \defeq \bG^D(\bx,\by^I)\cdot y_3(-f_1,-f_2,f_3)^T \\
\bu^{L1} &= -\frac{1}{2} \left( x_3\nabla_x-\hat{\bx}_3 \right) \phi^S(\bx), \quad\text{ with }\phi^S(\bx) \defeq G^S(\bx,\by)f_3 + G^S(\bx,\by^I)(-f_3),  \\
\bu^{L2} &= \frac{1}{2}\nabla_x \phi^{SZ}(\bx), \quad\text{ with }\phi^{SZ}(\bx)\defeq\left[G^S(\bx,\by)(f_3y_3) + G^S(\bx,\by^I)(-f_3y_3)\right].
\end{align}
\end{subequations}
$\bu^S$ denotes the Stokes kernel sum, $\phi^D$ denotes the Laplace dipole sum, and $\phi^S,\phi^{SZ}$ denote two Laplace monopole sums. The values and gradients of $\phi^D$, $\phi^S$ and $\phi^{SZ}$ are computed at the target point $\bx$.
It is straightforward to verify that $\bu(\bx)$ is equivalent to the original image system in Eq.~(\ref{eq:imageGVG}).
In this new image system, the Stokeslet sum and the two Laplace monopole sums are obviously neutral. The Laplace dipole sum is intrinsically neutral, because each dipole source is the asymptotic limit of zero distance between equal and opposite charges.
Therefore, each of the 4 kernel sums can be separately periodized. 
Therefore we claim this new image system is applicable for non-periodic, singly periodic and doubly periodic systems.

Eq.~(\ref{eq:imageNew}) obviously keeps the same computational complexity of the underlying summation methods, which is $\Ocal(N\log N)$ for FFT-based methods and $\Ocal(N)$ for KIFMM.
Eq.~(\ref{eq:imageNew}) is also close-to-optimal because although one more $G^S$ kernel sum is invoked compared to the image system given by Eq.~(\ref{eq:imageGVG}), this $G^S$ kernel sum is usually much faster than the sum for the Stokeslet $\bJ$, because the kernel $G^S$ is a scalar while the Stokeslet $\bJ$ is a $3\times3$ tensor.

\section{Numerical Results}
\label{sec:results}
In this section we present numerical results using the periodic KIFMM method developed in our previous work \citep{PBCFMM2018}.
It works by splitting the infinite periodic domain into a near field and a far field.
The near field is directly summed by KIFMM while the far field is added through a precomputed Multipole-To-Local (M2L) operator applied to the near field calculation results. A single parameter $p$ controls the accuracy and cost the KIFMM method, by placing $p$ equivalent points per edge, in total $6(p-1)^2+2$ equivalent point sources for each cubic box in the octree in the KIFMM algorithm. 
Using larger $p$ gives better accuracy but has higher computational cost. Approximately, $p=10$ gives single precision accuracy and $p=16$ gives double precision accuracy. 
The M2L operators for the kernels $\bJ$, $G^S$, and $\bG^D$ are constructed with the doubly periodic formulation derived by \citet{Lindbo_Tornberg_2011} and \citet{Tornberg_2015}.
The computer program is based on the high-performance package PVFMM developed by \citet{malhotra_pvfmm_2015}.
Since the gradient of Laplace dipole potential and Laplace monopole potential are required in Eq.~(\ref{eq:imageNew})b,c,d,
we modified the corresponding kernels to generate the potentials simultaneously with the gradients in the final Source-To-Target (S2T), Local-To-Target (L2T), and Multipole-To-Target (M2T) stages of KIFMM \citep{ying_kernel-independent_2004}.
%This eliminated the necessity of performing a separate KIFMM for the gradient kernels.

We further optimized the implementation by combining Eq.~(\ref{eq:imageNew}b,c) into one KIFMM because they share the same operator $\left(x_3 \nabla_x - \hat{\bx}_3\right)$ in the final assembling step. This optimization is achieved by representing the M2M, M2L, and L2L operations in the tree with Laplace monopole densities only. 

\subsection{Timing results}
\label{subsec:Timing}
The computation is timed on a 12-core 3.6GHz Intel Xeon workstation. 
$97^3$ point forces are placed at a set of random source points. Each force component is randomly generated from a uniform distribution in $[-0.5,0.5]$, and each coordinate component of the source points is randomly generated from a lognormal distribution with standard parameters $(0.2,0.5)$. The target points are chosen to be a set of $97^3$ Chebyshev quadrature points. The source and target points are both scaled and shifted to fill the half unit cube $[0,1)^2\times[0,0.5)$. The no-slip boundary condition is imposed on the plane $x_3=0$ by placing the image points in the other half unit cube $[0,1)^2\times(-0.5,0]$. 

With the image method, to evaluate the velocity at $T$ points due to $S$ point forces, the Stokes FMM evaluates from $2S$ source points (including image points) to $T$ target points.
The Laplace dipole FMM is from $S$ to $T$ points, and the two Laplace FMMs are both from $2S$ (including image points) to $T$ points. The timing results are shown in Table~\ref{tab:timing}. The Stokes FMM evaluates Eq.~(\ref{eq:imageNew}a) and the Monopole 2 FMM evaluates Eq.~(\ref{eq:imageNew}d). The Dipole \& Monopole 1 FMM evaluates both Eq.~(\ref{eq:imageNew}b,c) in a single operation, as explained in the optimization mentioned above. 

\begin{table}
	\centering
	\caption{\label{tab:timing}Timing (in seconds) results for $97^3$ target and source
		points with doubly periodic boundary condition. $\tau_{tree}$ is the time to construct the octree for FMM, $\tau_{\Ncal}$ is the time for near-field evaluations and  $\tau_{\Fcal}$ is the time for far-field evaluations. Stokes, Dipole, Monopole 1 and Monopole 2 correspond to $\bu^S$, $\bu^D$, $\bu^{L1}$, and $\bu^{L2}$ in Eq.~(\ref{eq:imageNew}), respectively. $p$ is the number of equivalent points per cubic box edge in KIFMM, controlling the accuracy and cost of KIFMM.
		%Sum refers to the final stage of assembling all pieces together in Eq.~(\ref{eq:imageNew}). 
	}
	\begin{tabular}{c|ccc|ccc|ccc}
		\hline
		&\multicolumn{3}{c}{Stokes}&\multicolumn{3}{c}{Dipole \& Monopole 1}&\multicolumn{3}{c}{Monopole 2} \\
		\hline
		$p$ &$\tau_{tree}$ &$\tau_{\Ncal}$ & $\tau_{\Fcal}$
		& $\tau_{tree}$ & $\tau_{\Ncal}$ & $\tau_{\Fcal}$ 
		& $\tau_{tree}$ & $\tau_{\Ncal}$ & $\tau_{\Fcal}$ 
		\\
		\hline
		$6$  &0.75 & 1.56 & 0.026 & 0.81 & 1.08 & 0.021 & 0.81 & 1.02 & 0.022 \\
		$8$  &0.70 & 2.03 & 0.051 & 0.75 & 1.43 & 0.041 & 0.69 & 1.37 & 0.042 \\
		$10$ &0.73 & 2.64 & 0.085 & 0.69 & 1.78 & 0.068 & 0.62 & 1.55 & 0.066\\
		$12$ &0.73 & 3.83 & 0.14 & 0.67 & 2.20 & 0.099 & 0.90 & 2.01 & 0.098 \\
		$14$ &0.74 & 5.14 & 0.19 & 0.63 & 2.59 & 0.14 & 0.91 & 2.38 & 0.15 \\
		$16$ &0.86 & 7.95 & 0.28 & 0.78 & 3.19 & 0.18 & 0.95 & 2.82 & 0.18 \\
		\hline
	\end{tabular}
\end{table}

\subsection{Accuracy results}
\label{subsec:Accuracy}
\begin{figure}
	\centering
\includegraphics[width=\linewidth]{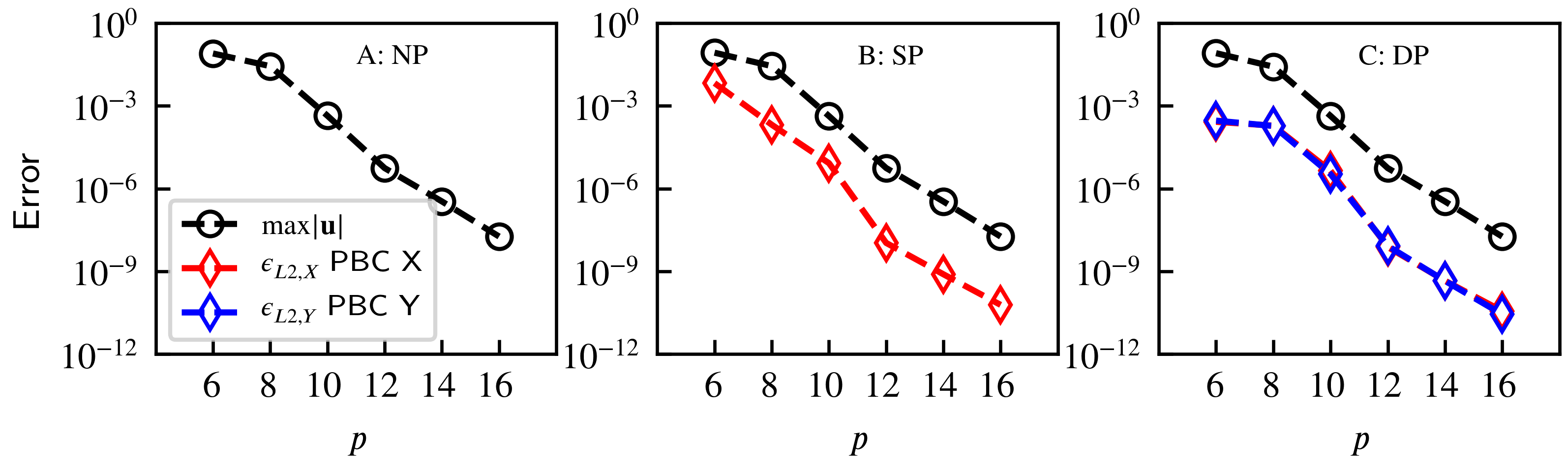}
\caption{\label{fig:accuracy} The accuracy test for the same set of source and target points specified in Section~\ref{subsec:Timing}. A: only no-slip condition is imposed on the $x_3=0$ plane. B: singly periodic boundary condition is imposed in $x_1$ direction. C: doubly periodic boundary condition is imposed in $x_1,x_2$ directions. In subplot C, $\epsilon_{L2,X}$ and $\epsilon_{L2,Y}$ almost overlap.}	
\end{figure}

In this section we use the same set of point forces and target points in the unit cube as in the last section, and present results for both non-periodic and doubly periodic boundary conditions.
To measure the accuracy of satisfying the no-slip boundary condition, a mesh of $97^2$ Chebyshev points is placed on the no-slip wall and the maximum component of flow velocity $\max\ABS{\bu}$ at these points is calculated.
To measure the accuracy of the doubly periodic boundary condition, a mesh of $97^2$ Chebyshev points is placed on each of the four side boundaries with $x_3\in[0,0.5)$.
%\footnote{PVFMM requires that the coordinates of source and target points cannot be equal to the exact floating point number $1.0$. Here we actually placed the points on the planes $x_1,x_2=1-10^{-12}$.}
The relative $L_2$ error for the flow velocity at the Chebyshev points on the two side walls $x_1=0,1,x_2\in[0,1),x_3\in[0,0.5)$ is used as the measurement of periodic boundary condition error in the $X$ direction, denoted by $\epsilon_{L2,X}$. 
Similarly $\epsilon_{L2,Y}$ is also measured. As shown in Figure.~\ref{fig:accuracy} the numerical errors in the no-slip and periodic boundary conditions converge exponentially with increasing $p$, and close to machine precision can be achieved. 

\section{The extension to other Stokes-related kernels}
\label{sec:extension}
\subsection{The Laplacian of Stokeslet.}
\label{subsec:StokesletLaplacian}
The Laplacian of the Stokeslet $\nabla_x^2\bJ$, is often invoked to compute the flow induced by a degenerate force doublet. For example, the flow induced by a slender fiber with length $L$ can be written as $\bu(\bx) = \int_{s=0}^{L} \left(1+\frac{1}{2}\epsilon^2\nabla_x^2\right)\bJ\left[\bx,\by(\bs)\right] \bff(s) $, where $\bff$ is the force density on the fiber, and $\epsilon\ll1$ is the slenderness parameter of the fiber \citep{Gotz2002,Nazockdast_Rahimian_Zorin_Shelley_2017}.
It is straightforward to construct the image system for $\nabla_x^2\bJ$ by directly taking the Laplacian of Eq.~(\ref{eq:imageNew}). 
However, a simpler image system can be derived by realizing that $\nabla_x^2\bJ(\bx,\by) = \nabla_x\bG^D(\bx,\by)$.
We denote this kernel by $\bQ(\bx,\by)$ because it is similar to the Green's function of an electrostatic quadrupole:
\begin{align}
	\bQ(\bx,\by) = \nabla_x^2\bJ(\bx,\by) = \nabla_x\bG^D(\bx,\by) 
%= \frac{1}{4\pi}\frac{1}{\ABS{\bx-\by}^5}\left(\ABS{\bx-\by}^2\bI - 3(\bx-\by)(\bx-\by)\right).
\end{align}
It is straightforward to verify that $\nabla_x^2\bQ = 0$ and $\nabla_x\cdot\bQ=0$.

Following the Papkovich-Neuber approach as used by \citet{Gimbutas_Greengard_Veerapaneni_2015}, we find a harmonic potential $\phi^{GD} \defeq \hat{\bx}_3\cdot[\bQ(\bx,\by^I)\bff^I]$ to construct the solution to Stokes equation to complete the image system. We choose the superscript $GD$ because $\bQ$ is the gradient of the Laplace dipole kernel $\bG^D$. 
%with $\bu=x_3\nabla \phi^Q-\hat{\bx}_3\phi^Q,p=2{\partial\phi^Q}/{\partial x_3}$ 
It is straightforward to verify that the following image system satisfies the no-slip boundary condition and the Stokes equation:
\begin{align}
\bu^{GD}(\bx) &= x_3 \nabla_x\phi^{GD}(\bx) - \hat{\bx}_3\cdot\phi^{GD}(\bx),\quad 
\phi^{GD}(\bx) \defeq \hat{\bx}_3 \left[\bQ(\bx,\by^I)\bff^I\right]\\
\bu &= \bQ(\bx,\by)\bff - \bQ(\bx,\by^I)\bff^I - 2\bu^{GD}, \quad p = -4\dpone{\phi^{GD}(\bx)}{x_3}
\end{align}
This image system can be calculated with two kernel sums, $\bQ(\bx,\by)\bff$ and $\bQ(\bx,\by^I)\bff^I$, both involving the $\bQ$ kernel only.
Further we do not need to rearrange and neutralize the sums, because $\bQ(\bx,\by)=\nabla_x \bG^D(\bx,\by)$ and $\bG^D$ is intrinsically neutral since it represents a Laplace dipole as discussed before. 
The doubly periodic formulation by \citet{Bleibel_2012} or the general formulation for electrostatic systems by \citet{Tornberg_2015} can be directly used as black-box periodic summation routines for the kernel $\bQ$.

\subsection{The Rotne-Prager-Yamakawa tensor}
\label{subsec:StokesletRPY}
The Rotne-Prager-Yamakawa (RPY) tensor \citep{Rotne_Prager_1969,Yamakawa_1970} is widely used in simulations of Brownian suspensions
and Langevin dynamics of biomolecules because it is a reasonably accurate approximation to hydrodynamics in Stokes flow and, more importantly, is designed to be always symmetric positive definite.
Without the wall, the RPY tensor is constructed from the Stokeslet $\bJ(\bx,\by)$ as:
\begin{align}
\label{eq:RPYfreespace}
\bu = \left(1+\frac{1}{6}a^2\nabla_x^2\right)\bu', \quad \bu' = \left(1+\frac{1}{6}b^2\nabla_y^2\right) \bJ(\bx,\by) \bff .
\end{align}
Here $\bu'$ is the velocity disturbance induced by the force $\bff$ on a `source' particle with radius $b$, and $\bu$ is the velocity of a `target' particle with radius $a$ induced by $\bu'$. The Laplacian terms $\frac{1}{6}a^2\nabla_x^2$ and $\frac{1}{6}b^2\nabla_y^2$ represent the finite-size effects of the target and source particles, respectively.
In the absence of the wall, the bi-Laplacian term $\nabla_x^2\nabla_y^2\bJ(\bx,\by)$ is zero.
However, the bi-Laplacian term is not zero in the presence of the wall, and the image system is significantly more complicated than Eq.~(\ref{eq:RPYfreespace}), as shown by \citet{Swan_Brady_2007}.
In this section we derive an image system for Eq.~(\ref{eq:RPYfreespace}) based on Eq.~(\ref{eq:imageNew}). This new image system is applicable with either periodic or non-periodic geometries. 

Starting from Eq.~(\ref{eq:imageNew}), the velocity $\bu$ of the target particle at $\bx$ induced by the force $\bff$ on the source particle located at $\by$ is:
\begin{align}
\bu = \left(1+\frac{1}{6}a^2\nabla_x^2\right)\bu', \quad \bu' = \left(1+\frac{1}{6}b^2\nabla_y^2\right)\left( \bu^S + \bu^D +\bu^{L1} +\bu^{L2} \right),
\end{align}
This can be simplified by realizing that $\nabla_y^2\bu^{L1}=\bm{0}$, $\nabla_x^2\bu^{L2}=\bm{0}$, and $\nabla_x^2\nabla_y^2\bu^S=\bm{0}$:
\begin{align}
\bu &= \bu^S + \bu^D +\bu^{L1} +\bu^{L2} + \frac{b^2}{6}\nabla_y^2\left(\bu^S + \bu^D + \bu^{L2}\right) \nonumber \\
& + \frac{a^2}{6}\nabla_x^2\left(\bu^S + \bu^D +\bu^{L1}\right) + \frac{a^2b^2}{36}\nabla_x^2\nabla_y^2\bu^D.
\end{align}
The $\nabla_x^2$ terms can be written as:
\begin{align}
	&\nabla_x^2 \bu^S= \bQ(\bx,\by)\bff_{xy} + \bQ(\bx,\by^I)(-\bff_{xy}), \\
	&\nabla_x^2 \bu^D= x_3 \nabla_x \nabla_x^2\phi^D + 2 \dpone{}{x_3} \nabla_x \phi^D-[0,0,\nabla_x^2\phi^D] = 2\dpone{}{x_3} \nabla_x \phi^D, \\
	&\nabla_x^2 \bu^{L1}= -\frac{1}{2}x_3 \nabla_x \nabla_x^2\phi^S - \dpone{}{x_3} \nabla_x \phi^S-[0,0,\nabla_x^2\phi^S] = -\dpone{}{x_3} \nabla_x \phi^S,
\end{align}
where we have utilized the fact that $\nabla_x^2\phi^D=0$, and $\nabla_x^2\phi^S=0$.

By symmetry $\nabla_y^2\bu^S =\nabla_x^2\bu^S$. Other $\nabla_y^2$ and the bi-Laplacian terms can be written as:
\begin{align}
	&\nabla_y^2\bu^{L2} = \nabla_x\phi^{DZ},\quad \phi^{DZ}\defeq\bG^D(\bx,\by)\cdot(0,0,f_3)^T + \bG^D(\bx,\by^I)\cdot(0,0,f_3)^T, \\
	&\nabla_y^2\bu^D = x_3\nabla_x\nabla_y^2\phi^D - (0,0,\nabla_y^2\phi^D),\\
	&\nabla_x^2\nabla_y^2 \bu^D = 2\dpone{}{x_3}\nabla_x \nabla_y^2\phi^D, \quad  \nabla_y^2\phi^D\defeq \phi^Q = \bG^Q(\bx,\by) \colon 2\begin{bmatrix}
	f_3 & 0 & 0 \\
	0 & f_3 & 0 \\
	f_1 & f_2 & 0 
	\end{bmatrix},
\end{align}
where $\nabla_y^2\phi^D$ is represented by the field $\phi^Q$ induced by a Laplace quadrupole with kernel function $\bG^Q(\bx,\by) = \nabla_y \bG^D(\bx,\by)$.
Here the symbol `$\colon$' denotes double-contraction of $3\times3$ tensors.
The image system for the RPY tensor in non-periodic and periodic geometries can be universally represented as a combination of the terms derived above. We summarize the results for monodisperse and polydisperse systems separately in the following, because in monodisperse systems the image system can be simplified with $a=b$.

\subsubsection{The kernel sums for monodisperse systems}
Using that in monodisperse systems $b=a$ for all particles, the particle radius $a$ is scaled out of the source strengths and 6 kernel sums are needed as shown in Table~\ref{tab:RPYmono}. 
\begin{table}
	\centering
	\caption{\label{tab:RPYmono} The kernel sums for the image RPY tensor in a monodisperse system. 
		%The velocity of the target particle can be assembled as Eq.~\ref{eq:FMMRPYmono}.
	}
	\begin{tabular}{l|c|c}
		\hline
		  & Source strength and location & Target values \\
		\hline
		Stokes 1 $\bu^S$ & $\bff_{xy}$ at $\by$, and $-\bff_{xy}$ at $\by^I$  & $\bu^S,\nabla_x^2\bu^S$ \\
		Laplace Monopole 1 $\phi^S$ & $f_3$ at $\by$, and $-f_3$ at $\by^I$  & $\phi^S, \nabla_x\phi^S, \nabla_x\nabla_x\phi^S$\\
		Laplace Monopole 2 $\phi^{SZ}$ & $f_3y_3$ at $\by$, and $-f_3y_3$ at $\by^I$   & $\phi^{SZ}, \nabla_x\phi^{SZ}, \nabla_x\nabla_x\phi^{SZ}$\\
		Laplace Dipole 1 $\phi^D$ & $y_3(-f_1,-f_2,f_3)^T$ at $\by^I$ &  $\phi^D, \nabla_x\phi^D, \nabla_x\nabla_x\phi^D$ \\
		Laplace Dipole 2 $\phi^{DZ}$ & $(0,0,f_3)^T$ at $\by$, and $(0,0,f_3)^T$ at $\by^I$ & $\nabla_x\phi^{DZ}$ \\
		Laplace Quadrupole $\phi^Q$ & 2$\begin{bmatrix}
		f_3 & 0 & 0 \\
		0 & f_3 & 0 \\
		f_1 & f_2 & 0 
		\end{bmatrix}$ at $\by$ & $\phi^Q$, $\nabla_x \phi^Q$ \\
		\hline
	\end{tabular}
\end{table}

The velocity of the target particle can be written as a combination of the target values:
\begin{align}\label{eq:FMMRPYmono}
\bu &= \left(1+\frac{a^2}{3}\nabla_x^2  \right) \bu^S + \left(x_3\nabla_x-\hat{\bx}_3+\frac{a^2}{3}\dpone{}{x_3}\nabla_x\right)\phi^D \nonumber -\frac{1}{2}\left(x_3\nabla_x-\hat{\bx}_3 + \frac{a^2}{3}\dpone{}{x_3}\nabla_x\right)\phi^S \nonumber \\ &+\frac{1}{2}\nabla_x\phi^{SZ}+\frac{a^2}{6}\nabla_x\phi^{DZ} + \left(\frac{a^2}{6} x_3\nabla_x-\frac{a^2}{6}\hat{\bx}_3+\frac{a^4}{18}\dpone{ }{x_3}\nabla_x\right) \phi^Q
\end{align}

\subsubsection{The kernel sums for polydisperse systems}
In this case we cannot simply scale the kernel sum results with the particle radius as in the previous case because the radius $b$ can be different for all source particles. Therefore $b$ must be included in the source strengths and 7 kernel sums are needed, as shown in Table~\ref{tab:RPYpoly}.

\begin{table}
	\centering
	\caption{\label{tab:RPYpoly}  The kernel sums for the image RPY tensor in a polydisperse system. 
	%	The velocity of the target particle can be assembled as Eq.~\ref{eq:FMMRPY}.
	}
	\begin{tabular}{l|c|c}
		\hline
		& Source strength and location  & Target values \\
		\hline
		Stokes 1 $\bu^S$ & $\bff_{xy}$ at $\by$, and $-\bff_{xy}$ at $\by^I$   & $\bu^S,\nabla_x^2\bu^S$ \\
		Stokes 2 $\bu_b^S$ & $b^2\bff_{xy}$ at $\by$, and $-b^2\bff_{xy}$ at $\by^I$   & $\nabla_x^2\bu_b^S$ \\
		Laplace Monopole 1 $\phi^S$ & $f_3$ at $\by$, and $-f_3$ at $\by^I$  & $\phi^S, \nabla_x\phi^S, \nabla_x\nabla_x\phi^S$\\
		Laplace Monopole 2 $\phi^{SZ}$ & $f_3y_3$ at $\by$, and $-f_3y_3$ at $\by^I$  & $\phi^{SZ}, \nabla_x\phi^{SZ}, \nabla_x\nabla_x\phi^{SZ}$\\
		Laplace Dipole 1 $\phi^D$ & $y_3(-f_1,-f_2,f_3)^T$ at $\by^I$  & $\phi^D, \nabla_x\phi^D, \nabla_x\nabla_x\phi^D$ \\
		Laplace Dipole 2 $\phi_b^{DZ}$ & $b^2(0,0,f_3)^T$ at $\by$, and $b^2(0,0,f_3)^T$ at $\by^I$ &  $\nabla_x\phi_b^{DZ}$ \\
		Laplace Quadrupole $\phi_b^Q$ & 2 $b^2\begin{bmatrix}
		f_3 & 0 & 0 \\
		0 & f_3 & 0 \\
		f_1 & f_2 & 0 
		\end{bmatrix}$ at $\by$ & $\phi_b^Q$, $\nabla_x \phi_b^Q$ \\
		\hline
	\end{tabular}
\end{table}

The velocity of the target particle is a combination of the target values:
\begin{align}\label{eq:FMMRPY}
	\bu &= \left(1+\frac{a^2}{6}\nabla_x^2 \right) \bu^S + \frac{1}{6}\nabla_x^2\bu_b^S + \left(x_3\nabla_x-\hat{\bx}_3+2\frac{a^2}{6}\dpone{}{x_3}\nabla_x\right)\phi^D -\frac{1}{2}\left(x_3\nabla_x-\hat{\bx}_3 +2\frac{a^2}{6}\dpone{}{x_3}\nabla_x\right)\phi^S \nonumber \\
	& +\frac{1}{2}\nabla_x\phi^{SZ}+\frac{1}{6}\nabla_x\phi_b^{DZ}+ \left(\frac{1}{6} x_3\nabla_x-\frac{1}{6}\hat{\bx}_3+\frac{a^2}{18}\dpone{ }{x_3}\nabla_x\right) \phi_b^Q
\end{align}

\begin{figure}
	\centering
	\includegraphics[width=\linewidth]{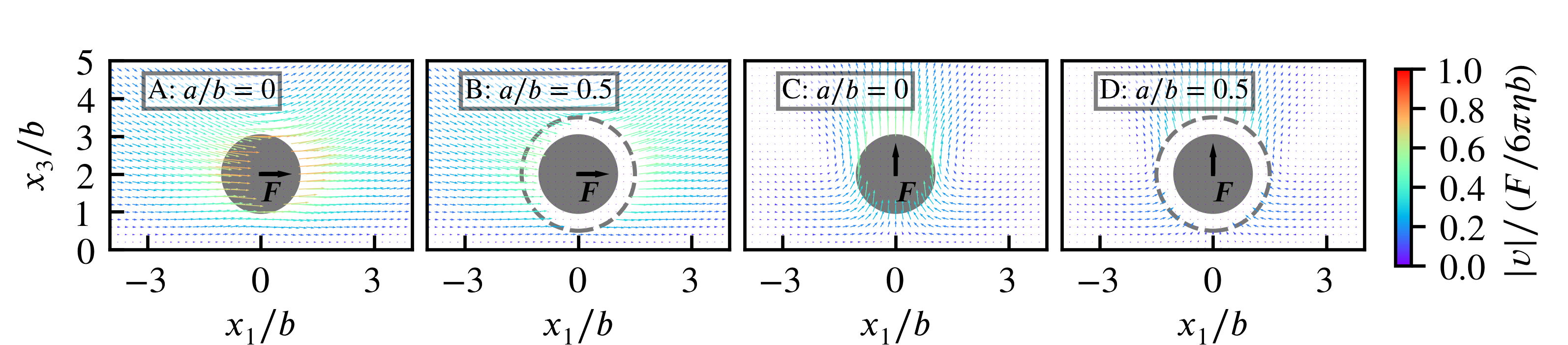}
	\caption{\label{fig:rpyflow} The velocity of a target particle induced by a source particle.
	The particle shown in the figure is the source particle with radius $b$, and a force $F$ is applied on the source particle in $\hat{\bx}_1$ (A,B) or $\hat{\bx}_3$ (C,D) directions. The vector field shows the velocity of the target particle with radius $a$ centered at some point in space. A and C: $a/b=0$, the target particle is infinitely small. B and D: $a/b=0.5$, the target particle center is not allowed to enter the dashed circle with radius $a+b$ due to non-overlap condition.}
\end{figure}

As we discussed in Section~\ref{sec:formulation}, the image systems tabulated in Table~\ref{tab:RPYmono} and~\ref{tab:RPYpoly} are universally applicable for non-periodic, singly periodic, and doubly periodic geometries, because the Stokes and Laplace monopole kernel sums are designed to be neutral, and the Laplace dipole and quadrupole kernel sums are intrinsically neutral.
Also the calculation of gradients and Laplacians ($\nabla_x\bu^S$, etc.) requires little extra cost in addition to the evaluation of field values ($\bu^S$, etc.) in both FMM and FFT type methods.
In FMM the gradients and Laplacians are straightforward to calculate in the final L2T, S2T and M2T stages because we can directly use the equivalent sources in KIFMM or the multipole basis functions in classic FMM to evaluate the gradients and Laplacians.
In FFT type methods the gradients and Laplacians can also be conveniently evaluated using interpolation schemes in the final stage where the kernel sum values are interpolated from the regular FFT mesh to the target points. The detailed cost analysis depends on the specific kernel sum methods used, but in any cases the summation scheme shown in Table~\ref{tab:RPYmono} and~\ref{tab:RPYpoly} does not change the algorithmic complexity of the underlying kernel sum methods. The final combination stages of Eq.~(\ref{eq:FMMRPYmono}) and~(\ref{eq:FMMRPY}) have negligible cost, similar to the case for the Stokeslet shown in Table~\ref{tab:timing}.

We conclude this section with a demonstration of the image system Eq.~(\ref{eq:FMMRPY}) in the polydisperse case, in the non-periodic geometry above a no-slip wall. The velocity of the target particle with radius $a$ for two cases in Fig.~\ref{fig:rpyflow}: $a/b=0$ in subplots A \& C and $a/b=0.5$ in subplots B \& D, are shown.
The no-slip wall is placed at the $x_3=0$ plane. The vector field shows the velocity of the target particle, generated by a force $\bF$ on the source particle $b$ located in $(0,0,2b)$.
The velocity vectors are colored by their magnitude. $\bF=(F,0,0)$ in A and B, and $\bF=(0,0,F)$ in C and D.
The dashed circles in B and D have radius $a+b$, showing the region where the center of the target particle cannot enter since we do not allow the source and target spheres to overlap. 
It is clear that the no-slip condition is satisfied.

\section{Discussion \& conclusion}
\label{conclusion}
In this paper we presented a new image system of a Stokeslet above a no-slip wall. 
This new image system generates exactly the same flow field in the non-periodic geometry as the image system derived by \citet{Blake_1971} and by \citet{Gimbutas_Greengard_Veerapaneni_2015}. Moreover, this new image system is straightforward to periodize with the partially periodic kernel summation methods, because the image system has been rearranged into neutral systems.
In other words, the periodic (or non-periodic) kernel sum methods can be called as black-box routines.
We demonstrated the accuracy and efficiency of this new image system in Section~\ref{sec:results}, using the periodic KIFMM by \citet{PBCFMM2018} as the black-box kernel summation routines.
Other summation methods can also be straightforwardly used without modifications \citep{Lindbo_Tornberg_2011,Lindbo_Tornberg_2011a,Nestler_Pippig_Potts_2015,Nestler_2016,barnett_unified_2016}.

The decomposition into neutral systems presented in Eq.~(\ref{eq:imageNew}) is not unique.
Following the idea of \citet{Tornberg_Greengard_2008}, the summations involving kernels $\bJ$ and $\bG^D$ can both be represented with several kernel sums involving the Laplace monopole kernel $G^S$ only.
This  reformulation allows the usage of classic Laplace FMM and Laplace Ewald methods.
The Stokeslet image system presented in this work can also be extended to the regularized Stokeslet, by directly integrating Eq.~(\ref{eq:imageNew}) over the regularization `blob' functions.
Then the regularized image system can be periodized, utilizing the doubly periodic schemes by \citet{Cortez_Hoffmann_2014} for the regularized Stokeslet,

The image systems for other fundamental solutions in Stokes flow have also been studied, including the doublet, stresslet, and rotlet kernels \citep{Gimbutas_Greengard_Veerapaneni_2015}. However they do not need to be rearranged into neutral systems because a doublet, a stresslet, or a rotlet is intrinsically a neutral force couple in Stokes flow.

In this work we also extended our image system of Stokeslet to its Laplacian and the RPY tensor.
In both cases the image systems is universally applicable for non-periodic, singly periodic or doubly periodic geometries.
The image system for $\nabla_x^2\bJ$ can be used in the simulations of rigid or flexible fibers above a wall in Stokes flow.
The image RPY tensor maintains the Symmetric-Positive-Definiteness of the mobility matrix, and can be widely used in Brownian dynamics simulations where the Brownian fluctuations must be generated from a SPD mobility matrix according to the fluctuation-dissipation theorem. 

\section{Acknowledgment}
Wen Yan thanks Aleksandar Donev for discussion of the RPY tensor.
MJS acknowledges support from NSF grants DMS-1463962 and DMS-1620331.

\section*{Appendix: Implementation}
The method mentioned in this paper has been implemented in the software \texttt{PeriodicFMM}, freely available on GitHub: \url{https://github.com/wenyan4work/PeriodicFMM}. The package is native in \texttt{C++} with interfaces in \texttt{C}, \texttt{Fortran}, and \texttt{Python}, and is fully parallelized with both \texttt{OpenMP} and \texttt{MPI}. The package is based on the author's fork of \texttt{PVFMM} \cite{PBCFMM2018}, also available on GitHub: \url{https://github.com/wenyan4work/pvfmm}.

\bibliographystyle{model1-num-names}

\begin{thebibliography}{33}
	\expandafter\ifx\csname natexlab\endcsname\relax\def\natexlab#1{#1}\fi
	\providecommand{\url}[1]{\texttt{#1}}
	\providecommand{\href}[2]{#2}
	\providecommand{\path}[1]{#1}
	\providecommand{\DOIprefix}{doi:}
	\providecommand{\ArXivprefix}{arXiv:}
	\providecommand{\URLprefix}{URL: }
	\providecommand{\Pubmedprefix}{pmid:}
	\providecommand{\doi}[1]{\href{http://dx.doi.org/#1}{\path{#1}}}
	\providecommand{\Pubmed}[1]{\href{pmid:#1}{\path{#1}}}
	\providecommand{\bibinfo}[2]{#2}
	\ifx\xfnm\relax \def\xfnm[#1]{\unskip,\space#1}\fi
	%Type = Article
	\bibitem[{Lele et~al.(2011)Lele, Swan, Brady, Wagner, and
		Furst}]{Lele_Swan_Brady_Wagner_Furst_2011}
	\bibinfo{author}{P.~P. Lele}, \bibinfo{author}{J.~W. Swan},
	\bibinfo{author}{J.~F. Brady}, \bibinfo{author}{N.~J. Wagner},
	\bibinfo{author}{E.~M. Furst},
	\newblock \bibinfo{title}{Colloidal diffusion and hydrodynamic screening near
		boundaries},
	\newblock \bibinfo{journal}{Soft Matter} \bibinfo{volume}{7}
	(\bibinfo{year}{2011}) \bibinfo{pages}{6844}.
	%Type = Article
	\bibitem[{Michailidou et~al.(2013)Michailidou, Swan, Brady, and
		Petekidis}]{Michailidou_Swan_Brady_Petekidis_2013}
	\bibinfo{author}{V.~N. Michailidou}, \bibinfo{author}{J.~W. Swan},
	\bibinfo{author}{J.~F. Brady}, \bibinfo{author}{G.~Petekidis},
	\newblock \bibinfo{title}{Anisotropic diffusion of concentrated hard-sphere
		colloids near a hard wall studied by evanescent wave dynamic light
		scattering},
	\newblock \bibinfo{journal}{The Journal of Chemical Physics}
	\bibinfo{volume}{139} (\bibinfo{year}{2013}) \bibinfo{pages}{164905}.
	%Type = Article
	\bibitem[{{Balboa Usabiaga} et~al.(2017){Balboa Usabiaga}, Delmotte, and
		Donev}]{Usabiaga_Delmotte_Donev_2017}
	\bibinfo{author}{F.~{Balboa Usabiaga}}, \bibinfo{author}{B.~Delmotte},
	\bibinfo{author}{A.~Donev},
	\newblock \bibinfo{title}{Brownian dynamics of confined suspensions of active
		microrollers},
	\newblock \bibinfo{journal}{The Journal of Chemical Physics}
	\bibinfo{volume}{146} (\bibinfo{year}{2017}) \bibinfo{pages}{134104}.
	%Type = Article
	\bibitem[{Kaya and Koser(2012)}]{Kaya_Koser_2012}
	\bibinfo{author}{T.~Kaya}, \bibinfo{author}{H.~Koser},
	\newblock \bibinfo{title}{Direct upstream motility in escherichia coli},
	\newblock \bibinfo{journal}{Biophysical Journal} \bibinfo{volume}{102}
	(\bibinfo{year}{2012}) \bibinfo{pages}{1514--1523}.
	%Type = Article
	\bibitem[{Spagnolie and Lauga(2012)}]{Spagnolie_Lauga_2012}
	\bibinfo{author}{S.~E. Spagnolie}, \bibinfo{author}{E.~Lauga},
	\newblock \bibinfo{title}{Hydrodynamics of self-propulsion near a boundary:
		predictions and accuracy of far-field approximations},
	\newblock \bibinfo{journal}{Journal of Fluid Mechanics} \bibinfo{volume}{700}
	(\bibinfo{year}{2012}) \bibinfo{pages}{105\u2013147}.
	%Type = Article
	\bibitem[{Wioland et~al.(2016)Wioland, Lushi, and
		Goldstein}]{Wioland_Lushi_Goldstein_2016}
	\bibinfo{author}{H.~Wioland}, \bibinfo{author}{E.~Lushi},
	\bibinfo{author}{R.~E. Goldstein},
	\newblock \bibinfo{title}{Directed collective motion of bacteria under channel
		confinement},
	\newblock \bibinfo{journal}{New Journal of Physics} \bibinfo{volume}{18}
	(\bibinfo{year}{2016}) \bibinfo{pages}{075002}.
	%Type = Article
	\bibitem[{Ezhilan and Saintillan(2015)}]{Ezhilan_Saintillan_2015}
	\bibinfo{author}{B.~Ezhilan}, \bibinfo{author}{D.~Saintillan},
	\newblock \bibinfo{title}{Transport of a dilute active suspension in
		pressure-driven channel flow},
	\newblock \bibinfo{journal}{Journal of Fluid Mechanics} \bibinfo{volume}{777}
	(\bibinfo{year}{2015}) \bibinfo{pages}{482–522}.
	%Type = Article
	\bibitem[{Blake(1971)}]{Blake_1971}
	\bibinfo{author}{J.~R. Blake},
	\newblock \bibinfo{title}{A note on the image system for a {Stokeslet} in a
		no-slip boundary},
	\newblock \bibinfo{journal}{Mathematical Proceedings of the Cambridge
		Philosophical Society} \bibinfo{volume}{70} (\bibinfo{year}{1971})
	\bibinfo{pages}{303--310}.
	%Type = Article
	\bibitem[{Mitchell and Spagnolie(2015)}]{Mitchell_Spagnolie_2015}
	\bibinfo{author}{W.~H. Mitchell}, \bibinfo{author}{S.~E. Spagnolie},
	\newblock \bibinfo{title}{Sedimentation of spheroidal bodies near walls in
		viscous fluids: glancing, reversing, tumbling and sliding},
	\newblock \bibinfo{journal}{Journal of Fluid Mechanics} \bibinfo{volume}{772}
	(\bibinfo{year}{2015}) \bibinfo{pages}{600\u2013629}.
	%Type = Article
	\bibitem[{Mitchell and Spagnolie(2017)}]{Mitchell_Spagnolie_2017}
	\bibinfo{author}{W.~H. Mitchell}, \bibinfo{author}{S.~E. Spagnolie},
	\newblock \bibinfo{title}{A generalized traction integral equation for stokes
		flow, with applications to near-wall particle mobility and viscous erosion},
	\newblock \bibinfo{journal}{Journal of Computational Physics}
	\bibinfo{volume}{333} (\bibinfo{year}{2017}) \bibinfo{pages}{462\u2013482}.
	%Type = Article
	\bibitem[{Gimbutas et~al.(2015)Gimbutas, Greengard, and
		Veerapaneni}]{Gimbutas_Greengard_Veerapaneni_2015}
	\bibinfo{author}{Z.~Gimbutas}, \bibinfo{author}{L.~Greengard},
	\bibinfo{author}{S.~Veerapaneni},
	\newblock \bibinfo{title}{Simple and efficient representations for the
		fundamental solutions of {Stokes} flow in a half-space},
	\newblock \bibinfo{journal}{Journal of Fluid Mechanics} \bibinfo{volume}{776}
	(\bibinfo{year}{2015}) \bibinfo{pages}{R1 (10 pages)}.
	%Type = Article
	\bibitem[{Nguyen and Leiderman(2015)}]{Nguyen_Leiderman_2015}
	\bibinfo{author}{H.-N. Nguyen}, \bibinfo{author}{K.~Leiderman},
	\newblock \bibinfo{title}{Computation of the singular and regularized image
		systems for doubly-periodic {Stokes} flow in the presence of a wall},
	\newblock \bibinfo{journal}{Journal of Computational Physics}
	\bibinfo{volume}{297} (\bibinfo{year}{2015}) \bibinfo{pages}{442--461}.
	%Type = Article
	\bibitem[{Hoffmann and Cortez(2017)}]{Hoffmann_Cortez_2017}
	\bibinfo{author}{F.~Hoffmann}, \bibinfo{author}{R.~Cortez},
	\newblock \bibinfo{title}{Numerical computation of doubly-periodic {Stokes}
		flow bounded by a plane with applications to nodal cilia},
	\newblock \bibinfo{journal}{Communications in Computational Physics}
	\bibinfo{volume}{22} (\bibinfo{year}{2017}) \bibinfo{pages}{620--642}.
	%Type = Article
	\bibitem[{Lindbo and Tornberg(2011)}]{Lindbo_Tornberg_2011}
	\bibinfo{author}{D.~Lindbo}, \bibinfo{author}{A.-K. Tornberg},
	\newblock \bibinfo{title}{Fast and spectrally accurate summation of 2-periodic
		{Stokes} potentials},
	\newblock \bibinfo{journal}{arXiv:1111.1815 [physics]}  (\bibinfo{year}{2011}).
	\bibinfo{note}{ArXiv: 1111.1815}.
	%Type = Article
	\bibitem[{Lindbo and Tornberg(2012)}]{Lindbo_Tornberg_2012}
	\bibinfo{author}{D.~Lindbo}, \bibinfo{author}{A.-K. Tornberg},
	\newblock \bibinfo{title}{Fast and spectrally accurate {Ewald} summation for
		2-periodic electrostatic systems},
	\newblock \bibinfo{journal}{The Journal of Chemical Physics}
	\bibinfo{volume}{136} (\bibinfo{year}{2012}) \bibinfo{pages}{164111}.
	%Type = Article
	\bibitem[{Shamshirgar and Tornberg(2017)}]{Shamshirgar_Tornberg_2017}
	\bibinfo{author}{D.~S. Shamshirgar}, \bibinfo{author}{A.-K. Tornberg},
	\newblock \bibinfo{title}{The spectral {Ewald} method for singly periodic
		domains},
	\newblock \bibinfo{journal}{Journal of Computational Physics}
	\bibinfo{volume}{347} (\bibinfo{year}{2017}) \bibinfo{pages}{341--366}.
	%Type = Article
	\bibitem[{Yan and Shelley(2018)}]{PBCFMM2018}
	\bibinfo{author}{W.~Yan}, \bibinfo{author}{M.~Shelley},
	\newblock \bibinfo{title}{Flexibly imposing periodicity in kernel independent
		fmm: A multipole-to-local operator approach},
	\newblock \bibinfo{journal}{Journal of Computational Physics}
	\bibinfo{volume}{355} (\bibinfo{year}{2018}) \bibinfo{pages}{214--232}.
	%Type = Article
	\bibitem[{Rotne and Prager(1969)}]{Rotne_Prager_1969}
	\bibinfo{author}{J.~Rotne}, \bibinfo{author}{S.~Prager},
	\newblock \bibinfo{title}{Variational treatment of hydrodynamic interaction in
		polymers},
	\newblock \bibinfo{journal}{The Journal of Chemical Physics}
	\bibinfo{volume}{50} (\bibinfo{year}{1969}) \bibinfo{pages}{4831--4837}.
	%Type = Article
	\bibitem[{Yamakawa(1970)}]{Yamakawa_1970}
	\bibinfo{author}{H.~Yamakawa},
	\newblock \bibinfo{title}{Transport properties of polymer chains in dilute
		solution{:} hydrodynamic interaction},
	\newblock \bibinfo{journal}{The Journal of Chemical Physics}
	\bibinfo{volume}{53} (\bibinfo{year}{1970}) \bibinfo{pages}{436--443}.
	%Type = Article
	\bibitem[{Hasimoto(1959)}]{Hasimoto_1959}
	\bibinfo{author}{H.~Hasimoto},
	\newblock \bibinfo{title}{On the periodic fundamental solutions of the {Stokes}
		equations and their application to viscous flow past a cubic array of
		spheres},
	\newblock \bibinfo{journal}{Journal of Fluid Mechanics} \bibinfo{volume}{5}
	(\bibinfo{year}{1959}) \bibinfo{pages}{317--328}.
	%Type = Article
	\bibitem[{Tornberg and Greengard(2008)}]{Tornberg_Greengard_2008}
	\bibinfo{author}{A.-K. Tornberg}, \bibinfo{author}{L.~Greengard},
	\newblock \bibinfo{title}{A fast multipole method for the three-dimensional
		{Stokes} equations},
	\newblock \bibinfo{journal}{Journal of Computational Physics}
	\bibinfo{volume}{227} (\bibinfo{year}{2008}) \bibinfo{pages}{1613--1619}.
	%Type = Article
	\bibitem[{Tornberg(2015)}]{Tornberg_2015}
	\bibinfo{author}{A.-K. Tornberg},
	\newblock \bibinfo{title}{The {Ewald} sums for singly, doubly and triply
		periodic electrostatic systems},
	\newblock \bibinfo{journal}{Advances in Computational Mathematics}
	\bibinfo{volume}{42} (\bibinfo{year}{2015}) \bibinfo{pages}{227--248}.
	%Type = Article
	\bibitem[{Malhotra and Biros(2015)}]{malhotra_pvfmm_2015}
	\bibinfo{author}{D.~Malhotra}, \bibinfo{author}{G.~Biros},
	\newblock \bibinfo{title}{{PVFMM}: {A} {Parallel} {Kernel} {Independent} {FMM}
		for {Particle} and {Volume} {Potentials}},
	\newblock \bibinfo{journal}{Communications in Computational Physics}
	\bibinfo{volume}{18} (\bibinfo{year}{2015}) \bibinfo{pages}{808--830}.
	%Type = Article
	\bibitem[{Ying et~al.(2004)Ying, Biros, and
		Zorin}]{ying_kernel-independent_2004}
	\bibinfo{author}{L.~Ying}, \bibinfo{author}{G.~Biros},
	\bibinfo{author}{D.~Zorin},
	\newblock \bibinfo{title}{A kernel-independent adaptive fast multipole
		algorithm in two and three dimensions},
	\newblock \bibinfo{journal}{Journal of Computational Physics}
	\bibinfo{volume}{196} (\bibinfo{year}{2004}) \bibinfo{pages}{591--626}.
	%Type = Article
	\bibitem[{G\"{o}tz(2002)}]{Gotz2002}
	\bibinfo{author}{T.~G\"{o}tz},
	\newblock \bibinfo{title}{On collocation schemes for integral equations arising
		in slender-body approximations of flow past particles with circular
		cross-section},
	\newblock \bibinfo{journal}{Journal of Engineering Mathematics}
	\bibinfo{volume}{42} (\bibinfo{year}{2002}) \bibinfo{pages}{203--221}.
	%Type = Article
	\bibitem[{Nazockdast et~al.(2017)Nazockdast, Rahimian, Zorin, and
		Shelley}]{Nazockdast_Rahimian_Zorin_Shelley_2017}
	\bibinfo{author}{E.~Nazockdast}, \bibinfo{author}{A.~Rahimian},
	\bibinfo{author}{D.~Zorin}, \bibinfo{author}{M.~Shelley},
	\newblock \bibinfo{title}{A fast platform for simulating semi-flexible fiber
		suspensions applied to cell mechanics},
	\newblock \bibinfo{journal}{Journal of Computational Physics}
	\bibinfo{volume}{329} (\bibinfo{year}{2017}) \bibinfo{pages}{173--209}.
	%Type = Article
	\bibitem[{Bleibel(2012)}]{Bleibel_2012}
	\bibinfo{author}{J.~Bleibel},
	\newblock \bibinfo{title}{{Ewald} sum for hydrodynamic interactions with
		periodicity in two dimensions},
	\newblock \bibinfo{journal}{Journal of Physics A: Mathematical and Theoretical}
	\bibinfo{volume}{45} (\bibinfo{year}{2012}) \bibinfo{pages}{225002}.
	%Type = Article
	\bibitem[{Swan and Brady(2007)}]{Swan_Brady_2007}
	\bibinfo{author}{J.~W. Swan}, \bibinfo{author}{J.~F. Brady},
	\newblock \bibinfo{title}{Simulation of hydrodynamically interacting particles
		near a no-slip boundary},
	\newblock \bibinfo{journal}{Physics of Fluids} \bibinfo{volume}{19}
	(\bibinfo{year}{2007}) \bibinfo{pages}{113306}.
	%Type = Article
	\bibitem[{Lindbo and Tornberg(2011)}]{Lindbo_Tornberg_2011a}
	\bibinfo{author}{D.~Lindbo}, \bibinfo{author}{A.-K. Tornberg},
	\newblock \bibinfo{title}{Spectral accuracy in fast ewald-based methods for
		particle simulations},
	\newblock \bibinfo{journal}{Journal of Computational Physics}
	\bibinfo{volume}{230} (\bibinfo{year}{2011}) \bibinfo{pages}{8744--8761}.
	%Type = Article
	\bibitem[{Nestler et~al.(2015)Nestler, Pippig, and
		Potts}]{Nestler_Pippig_Potts_2015}
	\bibinfo{author}{F.~Nestler}, \bibinfo{author}{M.~Pippig},
	\bibinfo{author}{D.~Potts},
	\newblock \bibinfo{title}{Fast {Ewald} summation based on {NFFT} with mixed
		periodicity},
	\newblock \bibinfo{journal}{Journal of Computational Physics}
	\bibinfo{volume}{285} (\bibinfo{year}{2015}) \bibinfo{pages}{280--315}.
	%Type = Article
	\bibitem[{Nestler(2016)}]{Nestler_2016}
	\bibinfo{author}{F.~Nestler},
	\newblock \bibinfo{title}{An {NFFT} based approach to the efficient computation
		of dipole-dipole interactions under various periodic boundary conditions},
	\newblock \bibinfo{journal}{Applied Numerical Mathematics}
	\bibinfo{volume}{105} (\bibinfo{year}{2016}) \bibinfo{pages}{25--46}.
	%Type = Article
	\bibitem[{Barnett et~al.(2016)Barnett, Marple, Veerapaneni, and
		Zhao}]{barnett_unified_2016}
	\bibinfo{author}{A.~H. Barnett}, \bibinfo{author}{G.~Marple},
	\bibinfo{author}{S.~Veerapaneni}, \bibinfo{author}{L.~Zhao},
	\newblock \bibinfo{title}{A unified integral equation scheme for
		doubly-periodic {Laplace} and {Stokes} boundary value problems in two
		dimensions},
	\newblock \bibinfo{journal}{arXiv:1611.08038 [math]}  (\bibinfo{year}{2016}).
	\bibinfo{note}{ArXiv: 1611.08038}.
	%Type = Article
	\bibitem[{Cortez and Hoffmann(2014)}]{Cortez_Hoffmann_2014}
	\bibinfo{author}{R.~Cortez}, \bibinfo{author}{F.~Hoffmann},
	\newblock \bibinfo{title}{A fast numerical method for computing doubly-periodic
		regularized {Stokes} flow in 3d},
	\newblock \bibinfo{journal}{Journal of Computational Physics}
	\bibinfo{volume}{258} (\bibinfo{year}{2014}) \bibinfo{pages}{1--14}.
	
\end{thebibliography}

\end{document}